\documentclass[12pt]{amsart}
\usepackage{geometry}                
\usepackage{graphicx}
\usepackage{amssymb}
\usepackage{epstopdf}

\usepackage{amsaddr}

\title{Construction of a universal ordinary differential equation $C^{\infty}$ of order 3}
\author{Etienne\ Couturier}
\address{Laboratoire\ "Matiere\ et\ Systemes\ Complexes"\ (MSC),\ UMR\ 7057\ CNRS,\ Universit\'e\ Paris\ 7\ Diderot,\ 75205\ Paris\ Cedex\ 13,\ France}
\author{Nicolas\ Jacquet}
\address{Lyc\'ee\ Jean-Baptiste\ Say,\ 11\ bis\ rue\ dÕAuteuil, 75016\ PARIS,\ France}
\begin{document}
\maketitle

\begin{abstract}
A universal ordinary differential equation $C^{\infty}$ of order 3 is constructed here. The equation is universal in the sense that any continuous function on a real segment can be approximated
by a solution of this equation with an arbitrary accuracy in uniform norm.
\end{abstract}



\section{Introduction and main result}
A most striking example of universal differential equation is due to Rubel (\cite{Rubel1981}). Inspired by results in theoretical computer science about universal Turing machine, he derived a universal algebraic differential equation of order 4. The algebraic formulation is very important because it means that the equation can be implemented using an analogic computer (\cite{Pour-El1974}). Rubel's method strongly relies on the singularity and the non-inversibility of the ADE which enables to glue many S-shaped solutions and thus to approximate any continuous function on a segment. His solutions are evidently not real-analytic functions. Using  similar methods, Elsner obtained a sequence $(P_n)$ of 3 th order ADE whose coefficients are not bounded. None of these equations is universal but any continuous function can be approximated with a given accuracy provided $n$ is sufficiently high. 
Universal ADE with real-analytic and even  polynomial solutions can be obtained at the cost of increasing the order of the equation (7th order for real-analytic solutions, 19th order  for polynomials)  (\cite{Boshernitzan1986}). \\
Despite all this impressive work very little attention has been paid to lipschitzian universal ODE; however many informations are already present in the existing litterature but are not clearly explicited. For instance, Boshernitzan has exhibited a four parameters  $(a,b,c,d)$ family of real function
 $$\forall x\in\mathbb{R},\ y(x)=c+\int_{a}^{x+a} \frac{bd}{1+d^2-\cos(bt)}\cos(\exp(t))dt $$
which are dense among the set of  continuous functions on a segment for the uniform norm. This universal formula can  easily be rewritten as a 4-dimensional system of $C^{\infty}$ ODE by turning two free parameters into constant variables as follows 

\begin{eqnarray*}
\dot{b}&=&0\\
\dot{d}&=&0\\
\dot{s}&=&1\\
\dot{y}&=&\frac{bd}{1+d^2-\cos(bs)}\cos(\exp(s))
\end{eqnarray*}
The two remaining parameters are included in the initial conditions $s(0)=a$ and $y(0)=c+\int_{0}^a \frac{bd}{1+d^2-\cos(bt)}\cos(\exp(t))dt $. Boshernitzan also provides a second formula which can similarly be rewritten as a system of 4 dimensional ODE.\\
As today system of lipshitzian ODE's are increasingly used as modeling tools especially in system biology: for instance a considerable literature is devoted to the robustness and the sensitivity of chemical network (\cite{Shinar2010}); it would be of very great use to clearly delimitate which ODE's are universal, which ones are close to universality and which ones are not. It would constitute a safeguard to the propensity for over-fitting and a robust ethical framework for modeling. High dimensional system of ODE's are often employed to model very tiny set of data leading to poorly robust estimation of parameters and initial conditions.\\
As a next step on this  project  a universal ODE which is $C^{\infty}$ of order 3 is constructed herein. The solutions are no more analytic-real as Boshernitzan's  but the order is lowered of one degree and the equation is of the simplest form $y'''=F(y,y',y'')$. \\

\noindent $\mathbf{Theorem.}$ For any $a<b$, there exists a third order  $C^{\infty}$  differential equation 
\begin{equation}
y'''=F(y,y',y'')
\label{Th}
\end{equation}
whose solutions are dense in $C^0([a,b])$ for the uniform norm.\\\\
Note the order 3 is the lower bound for a  lipschitzian  universal ODE.  Lipschitzian  ODE's of lower orders have too restricted space of solutions to be universal: first order lipschitzian equations $y'=F(y)$ only have monotonous solutions separated by the constant solutions which are the zeros of $F$; and second order lipschitzian equations $y''=F(y,y')$ have few types of solutions:  constant, monotonous, periodic and oscillating with a monotonous envelope. The second assertion can easily be shown knowing that the trajectories in the $(y,y')$ plane do not intersect.
 \\\\\\
 $\mathbf{Notations}$:\\ 
 - $Id$ stands for the function identity. \\
 - For $m\leqslant N$ two strictly positive integers, $ (i_1<i_2,\cdots< i_m) \in\ \{1,\cdots,N\}^m$ $m$ integers, $p_{(i_1,\cdots, i_m)}$ stands for the linear projection on the $m$ coordinates $(i_1,\cdots, i_m)$: 
 $$\forall x\in\ \mathbb{R}^N,\ p_{(i_1,\cdots, i_m)}(x_1,x_2,\cdots,x_N)=(x_{i_1},\cdots,x_{i_m})$$
- $dist$ stands for the euclidean distance in $\mathbb{R}^3$.\\
- For any $a<b$,  $f\in C^2([a,b])$
\begin{eqnarray*}
\forall t\in [a,b],\ \|f(t)\|&=&\sqrt{(f(t))^2+(f'(t))^2+(f''(t))^2}\\
\|f\|_s&=&\sup_{[a,b]}\|f\|\\
\|f\|_i&=&\inf_{[a,b]}\|f\|
\end{eqnarray*}
- For any $a<b$,  $f\in C^0([a,b])$, the support of $f$ is defined by:
$$Supp (f)=\{t\in [a,b],\ f(t)\neq0\}$$
- For any $a<b$, $k\geqslant 1$,  $f\in C^{k-1}([a,b])$, we note $f^{*k}=(f,f',\cdots,f^{(k-1)})$ and for any $j\geqslant k$,$f^{j*k}=(f^{j},f^{(j+1)},\cdots,f^{(k-1)})$\\
- For any $(x_0,y_0,z_0)\in \mathbb{R}^3$ and $\epsilon>0$, we note: $$\mathbf{B}((x_0,y_0,z_0),\epsilon)=\{(x,y,z)\in \mathbb{R}^3,\ \sqrt{(x-x_0)^2+(y-y_0)^2+(z-z_0)^2}<\epsilon\}$$
- We use the definition of  cite{Berger1992}. Let $a<b$, $\Gamma: [a,b]\to \mathbb{R}^3$ be a $C^1$ immersion. For $t\in [a,b]$, the normal space to $([a,b],\Gamma)$ at $t$ is defined by:
$$N_{\Gamma,t}[a,b]=\{z\in\ \mathbb{R}^3: (z|\Gamma'(t))=0\}$$
with $ (\ |\ )$ the standard scalar product of $\mathbb{R}^3$.\\
We set
\begin{eqnarray*}
N_\Gamma[a,b]&=&\{(t,v)\in\ [a,b]\times\mathbb{R}^3:\ v\in N_{\Gamma,t}[a,b]\}\\
N_\Gamma^\epsilon [a,b]&=&\{(t,v)\in\ N_\Gamma[a,b]:\ \sqrt{v_1^2+v_2^2+v_3^2}<\epsilon\}
\end{eqnarray*}
$N_\Gamma[a,b]$ is the normal fiber bundle to $([a,b],\Gamma)$.\\
 The canonical map $Can: N_\Gamma[a,b]\to \mathbb{R}^3$ is defined by:
\begin{eqnarray}
Can(t,v)=\Gamma(t)+v\label{Def_Can}
 \end{eqnarray}
The tubular neighborhood is defined by:
\begin{eqnarray}
Tub^\epsilon_\Gamma [a,b]=Can(N^\epsilon_\Gamma [a,b]).\label{Def_Tube}
 \end{eqnarray}
- For $a<b$, two real numbers $\gamma$,  for $\sigma\in\{0,1\}$, $t\in]a,b[$, and \\
$0<\beta<\min(b-t,t-a)$, we define $\psi_{\beta,t,\gamma}$ on $[a,b]$:
\begin{eqnarray*}
\forall u\in ]t-\beta,t+\beta[,&\:& \psi_{\beta,t,\gamma}(u)= \exp{\Big(-\frac{\gamma}{2}\frac{\beta^2(u-t)^2}{\beta^2-(u-t)^2}}\Big)\\
\forall u\in [a,t-\beta]\cup[t+\beta,b],&\:& \psi_{\beta,t,\gamma}(u)=0 
\end{eqnarray*}
A third order development of  $\psi_{\beta,t,\gamma}$ for $u$ in a neighborhood of $t$ gives:
\begin{eqnarray*}
\psi_{\beta,t,\gamma}(u)= 1-\frac{\gamma}{2}(u-t)^2+O((u-t)^4)\\
\end{eqnarray*}
\begin{eqnarray*}
\psi_{\beta,t,\gamma}^{*3}(t)=\left(
\begin{array}{ccc}
1& 0&-\gamma
\end{array}
\right)
\end{eqnarray*}
and
\begin{eqnarray*}
\psi_{\beta,t,\gamma}^{1*4}(t)
=\left(
\begin{array}{ccc}
0&-\gamma&0
\end{array}
\right)
\end{eqnarray*}
-  For $\epsilon>0$, we write:
\begin{eqnarray}
 \forall t\in[a,b],\ c_\epsilon(t)=\epsilon\cos(\frac{t}{\epsilon^2})
 \label{def_c_eps} 
 \end{eqnarray}\\\\

The main theorem will be proved using the following lemma $1,\ 2,\ 3$. The proofs of these lemma themselves are detailed in the next section and involve the lemma $4,\ 5,\ 6,\ 7,\ 8$ included in the appendix.\\\\

\noindent $\mathbf{Lemma\ 1.}$ For any $a<b$, there exists a sequence of $C^{\infty}$  function $(f_n)_{n\in\mathbb{N}}$ dense in $C^0([a,b])$ for the uniform norm such as for any integer $n\geqslant 0$
\begin{eqnarray}
\|f_n\|_s+1<\|f_{n+1}\|_i
\label{Hyp1_Lemme1}
\end{eqnarray}
\begin{eqnarray}
\forall t\in[a,b], (f_n'(t), f_n''(t)) \neq (0,0).
\label{Hyp2_Lemme1}
\end{eqnarray}

\noindent $\mathbf{Lemma\ 2.}$ For any $a<b$, there exists a sequence of $C^{\infty}$  function $(g_n)_{n\in\mathbb{N}}$ dense in $C^0([a,b])$ for the uniform norm  such as for any integer $n\geqslant 0$
\begin{eqnarray}
\|g_n\|_s+1<\|g_{n+1}\|_i
\label{Hyp1_Lemme2}
\end{eqnarray}
\begin{eqnarray}
\forall t\in[a,b],Ê (g_n'(t), g_n''(t)) \neq (0,0)
\label{Hyp2_Lemme2}
\end{eqnarray}
\begin{equation}
\forall\ c,d \in  [a,b] , c\neq d \Rightarrow g_n^{*3}(c)\neq g_n^{*3}(d).
\label{Hyp3_Lemme2}
\end{equation}
\\\\
\noindent $\mathbf{Lemma\ 3.}$ Let $g$ be a  $C^{\infty}$ function on $[a,b]$ such as
\begin{eqnarray}
\forall t\in[a,b],Ê (g'(t), g''(t)) \neq (0,0)
\label{Hyp1_Lemme3}
\end{eqnarray}
\begin{equation}
\forall\ c,d \in  [a,b] , c\neq d \Rightarrow g^{*3}(c)\neq g^{*3}(d).
\label{Hyp2_Lemme3}
\end{equation}
Then for any $\epsilon$ small enough, there exists a $C^{\infty}$ function $F$ of three variables whose support is included in $Tub^\epsilon_{g^{*3}}[a,b]$ (See (\ref{Def_Tube}) for a definition), such that
\begin{equation}
g'''=F(g,g',g'').
\label{Resultat_Lemme_1}
\end{equation}
\\\\
$\mathbf{Proof\ of\ the\ theorem.}$  From Lemma 2, we get a sequence  of functions $(g_n)$ dense in $C^0([a,b])$ verifying (\ref{Hyp1_Lemme2}, \ref{Hyp2_Lemme2},  \ref{Hyp3_Lemme2}). \\
For any integer $n$, $(g_n,g_n',g_n''):[a,b]\to \mathbb{R}^3$ is an immersion since for any $t\in\ [a,b],\ (g_n'(t),g_n''(t),g_n'''(t))\neq(0,0,0)$ by (\ref{Hyp2_Lemme2}); therefore the tubular neighborhood $Tub^\frac{1}{3}_{g_n^{*3}}[a,b]$ are well-defined for any $n$.\\
Moreover $(\ref{Hyp1_Lemme2})$ implies:
\begin{eqnarray*}
\forall t\in[a,b],\ \forall n\in  \mathbb{N}, \|g_{n}\|_i\leqslant\|g_{n}(t)\|&<&\|g_n\|_s+1<\|g_{n+1}\|_i\leqslant\|g_{n+1}(t)\|
\end{eqnarray*}
which implies by a trivial induction:
\begin{eqnarray*}
\forall t\in[a,b],\ \forall n\in  \mathbb{N},\forall m\in  \mathbb{N},n<m  \Rightarrow1\leqslant\|g_{m}(t)\|-\|g_n(t)\|\leqslant\|g_{m}(t)-g_n(t)\|
\end{eqnarray*} 
 \normalsize{These} inequalities imply that  the sequence of tubular neighborhoods  $Tub_{g_n^{*3}}^\frac{1}{3}[a,b]$ are all disjoint.\\
Using Lemma 3, we obtain a sequence of $C^{\infty}$ functions  $(F_n)$ such as for each $n>0$ the support of $F_n$ is included in $Tub_{g_n^{*3}}^\frac{1}{3}[a,b]$ and $F_n$ satisfies (\ref{eqn})
\begin{equation}
g_n'''=F_n(g_n,g_n',g_n'').
\label{eqn}
\end{equation}
We note
\begin{equation}
F=\sum_n F_n.
\end{equation}
As the supports of the $F_n$ are disjoint, the sum is locally finite and $F$ is  $C^{\infty}$ too.
By construction, for any $n\in \mathbb{N}$
\begin{equation}
g_n'''=F(g_n,g_n',g_n'').
\end{equation}
\\\\
\section{Proof of the Lemmas}
\noindent $\mathbf{Proof\ of\ Lemma\ 1.}$\\\\
Let $(P_n)$ be a sequence of non-constant polynomials dense in  $C^{0}([a, b])$ for the uniform norm. Let choose a sequence $(\epsilon_n)$ decreasing to $0$ such as for $n=1$
\begin{equation}
1+\|P_{1}\|_s
<\frac{1}{\epsilon_1} \\
\end{equation}
and for $n>1$
\begin{equation}
1+\|P_{n-1}\|_s+\sqrt{\frac{1}{\epsilon_{n-1}^6}+\epsilon_{n-1}^2}+\|P_{n}\|_s<\frac{1}{\epsilon_{n}}.
\label{eq:2}
\end{equation}
\\\\
 For an integer $n>0$, the sequence $(f_n)$ is defined by:
\begin{equation}
f_n=P_n+c_{\epsilon_n} \label{dodouii}
\end{equation}
See (\ref{def_c_eps}) for the definition of $c_{\epsilon_n}$.\\
As  for any $n$, $\epsilon_n<1$, $\frac{1}{\epsilon_n^2}<\frac{1}{\epsilon_n^6}$, $c_{\epsilon_n}$ is bounded as follows:\\\\
\small{\begin{eqnarray}
\|c_{\epsilon_n}\|_i&=\inf_{t\in\ [a,b]}\sqrt{(\epsilon_n^2+\frac{1}{\epsilon_n^6})\cos(\frac{t}{\epsilon_n^2})^2+\frac{1}{\epsilon_n^2}\sin(\frac{t}{\epsilon_n^2})^2}\geqslant&\frac{1}{\epsilon_n} \label{bound_inf}\\
\|c_{\epsilon_n}\|_s&=\sup_{t\in\ [a,b]}\sqrt{(\epsilon_n^2+\frac{1}{\epsilon_n^6})\cos(\frac{t}{\epsilon_n^2})^2+\frac{1}{\epsilon_n^2}\sin(\frac{t}{\epsilon_n^2})^2}\leqslant&\sqrt{\epsilon_n^2+\frac{1}{\epsilon_n^6}} \label{bound_sup}
\end{eqnarray}}

\noindent The triangular inequality  for $f_n$ reads:
\begin{eqnarray}
\|c_{\epsilon_n}\|_i-\|P_{n}\|_s\leqslant\|f_{n}\|\leqslant\|c_{\epsilon_n}\|_s+\|P_{n}\|_s \label{opa0}
\end{eqnarray}
Combining (\ref{opa0}) with the two bounds (\ref{bound_inf}, \ref{bound_sup}), it yields 
\begin{eqnarray}
\frac{1}{\epsilon_n}-\|P_{n}\|_s\leqslant\|f_{n}\|\leqslant\sqrt{\epsilon_n^2+\frac{1}{\epsilon_n^6}}+\|P_{n}\|_s \label{opa}
\end{eqnarray}
(\ref{eq:2}) can be rewritten:
\begin{eqnarray}
1+\sqrt{\epsilon_{n-1}^2+\frac{1}{\epsilon_{n-1}^6}}+\|P_{n-1}\|_s< \frac{1}{\epsilon_n}-\|P_{n}\|_s \label{dodo}
\end{eqnarray}
Using the inequality (\ref{opa}) expressed in $n-1$ on the right of (\ref{dodo}), and expressed in $n$ on the left of (\ref{dodo}) we get:
\begin{eqnarray*}
\|f_{n-1}\|_s+1<\|f_{n}\|_i 
\end{eqnarray*}
so that (\ref{Hyp1_Lemme1}) is verified by $f_n$.\\
The triangular inequalities  applied to $\sqrt{f_{n}^{\prime2}+f_{n}^{\prime\prime2}}$ combined to a trivial minoration of $\|P_{n}\|_s$ gives:
\begin{equation}
\frac{1}{\epsilon_n}-\|P_{n}\|_s\leqslant\frac{1}{\epsilon_n}-\sqrt{P_{n}^{\prime2}+P_{n}^{\prime\prime2}}\leqslant\sqrt{f_{n}^{\prime2}+f_{n}^{\prime\prime2}}
\label{ineq0}
\end{equation}
As the sequence $(P_n)$ is dense in $C^0([a,b])$ for the uniform norm, and $(\epsilon_n)$ is decreasing to $0$ at the infinity, the sequence $(f_n)$ is also dense in $C^0([a,b])$  for the uniform norm.\\
Combining (\ref{ineq0})  with (\ref{dodo}) yields to the following formula useful in the next proof:
\begin{equation}
1<\sqrt{f_{n}^{\prime2}+f_{n}^{\prime\prime2}}\label{ineq}
\end{equation}
and (\ref{Hyp2_Lemme1}) is satisfied by $f_n$.
\\\\
\noindent $\mathbf{Proof\ of\ Lemma\ 2.}$ Lemma 1 ensures the existence of a sequence $(f_n)_{n\in \mathbb{N}}$ which verifies (\ref{Hyp1_Lemme1}), (\ref{Hyp2_Lemme1}) such as for each $n$, $f_n$ is a sum of a non-constant polynomial $P_n$ and $c_{\epsilon_n}$. \\
\noindent We construct by induction two related sequences of function $(\delta_n)$  and $(g_n)$ linked by $g_n=f_n+\delta_n$; such as for $n\geqslant 0$, $g_n$ is verifying  both (\ref{Hyp2_Lemme2}) and (\ref{Hyp3_Lemme2}), while $\delta_n$ is verifying
\begin{eqnarray}
&&\|\delta_n\|_s<\frac{1}{n+1};\label{eququi1}
\end{eqnarray}
and such as for $n>0$, they also verified (\ref{Hyp1_Lemme2}) and
\begin{eqnarray}
&&\|f_{n}\|_i-\|g_{n-1}\|_s-1>0.\label{eququi2}
\end{eqnarray}
Let $n=0$, and $\xi_0=\|f_{1}\|_i-\|f_{0}\|_s-1$; $\xi_0>0$ by (\ref{Hyp1_Lemme1}). As $f_0=P_0+c_{\epsilon_0}$ and (\ref{Hyp2_Lemme1}) is verified, Lemma 8 ensures there exists a linear combination $\delta_0$ of plateau functions real-analytic on the interior of their support such as:
\begin{eqnarray}
\|\delta_0\|_s<\min(\frac{\xi_0}{10},1)
\label{conc2_Lemme9}
\end{eqnarray}
and such  that $g_0=f_0+\delta_0$ verifies both (\ref{Hyp2_Lemme2}), (\ref{Hyp3_Lemme2}). For $n=0$, (\ref{Hyp2_Lemme2}), (\ref{Hyp3_Lemme2}) and (\ref{eququi1}) are verified.
Moreover, (\ref{eququi2}) is verified for $n=1$:
\begin{eqnarray*}
\|f_{1}\|_i-\|g_{0}\|_s-1\geqslant \|f_{1}\|_i-\|f_{0}\|_s-\|\delta_{0}\|_s-1\geqslant \frac{9\xi_0}{10}>0
\end{eqnarray*}
\\
Let $n>0$ and suppose  a sequence $(g_{i})$ has been constructed for $i<n$ such as (\ref{Hyp2_Lemme2}), (\ref{Hyp3_Lemme2}) and (\ref{eququi1}) are verified and such as for any $i\in\{1,\cdots,n\}$, (\ref{Hyp1_Lemme2}) and (\ref{eququi2}) are also verified. Let $$\xi_n=\min(\|f_{n}\|_i-\|g_{n-1}\|_s-1,\|f_{n+1}\|_i-\|f_{n}\|_s-1).$$
$\xi_n>0$ by (\ref{Hyp1_Lemme1}) and (\ref{eququi2}). \\
As $f_n=P_n+c_{\epsilon_n}$ and (\ref{Hyp2_Lemme1}) is verified, Lemma 8 ensures there exists a linear combination $\delta_n$ of plateau functions real-analytic on their open support such as:
\begin{eqnarray}
\|\delta_n\|_s<\min(\frac{\xi_n}{10},\frac{1}{n+1});
\label{conc2_Lemme9}
\end{eqnarray}
and such  that $g_n=f_n+\delta_n$ verifies both (\ref{Hyp2_Lemme2}), (\ref{Hyp3_Lemme2}). The lower bound of $\|g_{n}\|_i$  obtained by triangular inequality gives:
$$\|g_{n}\|_i-\|g_{n-1}\|_s-1\geqslant (\|f_{n}\|_i-\|\delta_n\|_s)-\|g_{n-1}\|_s-1\geqslant \frac{9\xi_n}{10}>0$$
which means (\ref{Hyp1_Lemme2}) is also preserved for $n$.
Moreover the upper bound of $\|g_n\|_s$ obtained by triangular inequality gives:
 \begin{eqnarray*}
 \|f_{n+1}\|_i-\|g_n\|_s-1\geqslant \|f_{n+1}\|_i-(\|f_n\|_s+\|\delta_n\|_s)-1\geqslant \|f_{n+1}\|_i-\|f_n\|_s-\frac{\xi_n}{10} -1
 \end{eqnarray*}
 by definition of $\xi_n$
 \begin{eqnarray*}
 \|f_{n+1}\|_i-\|g_n\|_s-1\geqslant \frac{9 (\|f_{n+1}\|_i-\|f_n\|_s-1)}{10}>0
 \end{eqnarray*}
which proves the induction step for $(\ref{eququi2})$.\\
As both initial step and induction step have been proven, the two related sequences of function $(\delta_n)$  and $(g_n)$ verifying (\ref{Hyp1_Lemme2}), (\ref{Hyp2_Lemme2}), (\ref{Hyp3_Lemme2}) and (\ref{eququi1}) can be constructed for any $n>0$. As the sequence $(\|\delta_n\|)$ has for limit $0$ when $n$ approaches the infinity because of (\ref{eququi1}) and as the sequence $(f_n)$ is dense in $C^0([a,b])$, the sequence $(g_n)$ is also dense in $C^0([a,b])$ for the uniform norm. \\

\noindent $\mathbf{Proof\ of\ Lemma\ 3.}$ Hypothesis (\ref{Hyp1_Lemme2},\ref{Hyp2_Lemme2}) means that $g^{*3}=(g,g',g'')$ is a $C^1$ embedding from $[a,b]$ to $\mathbb{R}^3$. According to the theorem 2.7.12 of cite{Berger1992}  for $\epsilon$ sufficiently small, the canonical application $Can:N^\epsilon_{g^{*3}}[a,b] \to Tub^\epsilon_{g^{*3}}[a,b]$ is a $C^{1}$ diffeomorphism (See (\ref{Def_Can},\ref{Def_Tube}) for definitions).

\noindent Let's  define $F$ by its restriction.
\begin{eqnarray*}
 \forall (x_1,x_2,x_3)&\in&\ Tub^\epsilon_{g^{*3}}[a,b],\\
 F(x_1,x_2,x_3)&=&(g'''\circ p_1 \circ Can^{-1}(x_1,x_2,x_3)) \times  \exp\Big(\frac{1}{\epsilon^2}-\frac{1}{\epsilon^2-(p_{(2,3,4)} \circ Can^{-1}(x_1,x_2,x_3))^2}\Big)
\end{eqnarray*}
and 0 outside.\\
$F$ is $C^{\infty}$ by construction and $g$ satisfies the equation
\begin{equation}
g'''=F(g,g',g'')
\end{equation}
$\mathbf{Acknowledgment}$
Etienne Couturier is grateful both to Jean-Pierre Fran\c coise for introducing him to this subject, to Vincent Bansaye for a great help in the redaction. \\\\
\\
\section{Appendix}

\noindent $\mathbf{Lemma\ 4.}$ Let $a<b$, $\epsilon>0$, $P$ be a non-constant polynomial on $[a,b]$ and $$f=P+c_\epsilon.$$
\noindent We suppose:
\begin{eqnarray}
\forall t \in [a,b], (f'(t),f''(t))\neq (0,0)
\label{Hypo_lemme_4}
\end{eqnarray}
Let $l$ be the number of $(E_i)_{1\leqslant i\leqslant l} \in\mathbb{R}^3$ such that the equation for the variable $t$:
$$ E_i=f^{*3}(t)$$
admits $m_i$ solutions $(t_{i,j})_{1\leqslant j\leqslant m_i} \in [a,b]$ with $m_i>1$. Then $l$ is finite as well as the total number of solutions: $\Sigma_{i=1}^l m_i$. \\\\
\noindent $\mathbf{Proof\ of\ Lemma\ 4.}$ Let us suppose $\Sigma_{i=1}^l m_i$  to be non-finite. As the $t_{i,j}$ are included in a segment, there exists an accumulation point $(t'_{\infty,1},t'_{\infty,2})\in\ [a,b]^2$. Let  $(t'_{m,1},t'_{m,2})_{m\in \mathbb{N}}\in\ [a,b]^\mathbb{N}$ be  a injective sequence in $(t'_{m,1})$ and $(t'_{m,2})$ converging to $(t'_{\infty,1},t'_{\infty,2})$  such that: $$\forall m\in \mathbb{N},\  t'_{m,1}\neq t'_{m,2},\ f^{*3}(t'_{m,1})=f^{*3}(t'_{m,2})$$
We note:
\begin{eqnarray*}
F_0&:&(t_1,t_2)\in[a,b]^2\to (f(t_1)-f(t_2))\\ 
F_1&:&(t_1,t_2)\in[a,b]^2\to (f'(t_1)-f'(t_2))
\end{eqnarray*}

\noindent (\ref{Hypo_lemme_4}) ensures either $f'(t'_{\infty,2})\neq 0$ or  $f''(t'_{\infty,2})\neq 0$, let first suppose $f'(t'_{\infty,2})\neq 0$. \\
$f,\ f'$ are the sum of two real-analytic functions of one-variable so $F_0,\ F_1$  are real analytic functions of two variables. As $\frac{\partial F_0}{\partial t_2}(t'_{\infty,1},t'_{\infty,2})=-f'(t'_{\infty,2})\neq 0$, the implicit function theorem for real-analytic function (Theorem 2.3.1 of cite{Krantz1992}) ensures there exists a connected neighborhood $U$  of $t'_{\infty,1}$, a connected  neighborhood $V$ of $t'_{\infty,2}$, and a unique real analytic function $h:U\to V$  such that:
\begin{eqnarray}
\forall t \in\ U, \ F_0(t,h(t))=0
\label{F0}
\end{eqnarray}
\noindent Provided $m$ is sufficiently high $t'_{m,2}$ lays in $V$,  the unicity of $h$ ensures $$t'_{m,2}=h(t'_{m,1}).$$ As $(t'_{m,1})_{m\in \mathbb{N}}$ is an injective sequence, the function $t\to F_1(t,h(t))$ admits an accumulation of zeros around $t'_{\infty,1}$ so by analyticity: 
\begin{eqnarray}
\forall t\in U,\ F_1(t,h(t))=f'(t)-f'(h(t))=0
\label{F1_1}
\end{eqnarray}
Had we suppose $f''(t'_{\infty,2})\neq 0$ rather than $f'(t'_{\infty,2})\neq 0$,  the implicit function theorem would have been applied to $F_1$ rather than $F_0$ leading to the same output: 
\begin{eqnarray}
\forall t\in U,\ F_0(t,h(t))=F_1(t,h(t))=0 \label{oukoui}
\end{eqnarray}
Deriving (\ref{F0}) gives:
\begin{eqnarray}
\forall t\in U,\ \frac{dF_0(t,h(t))}{dt}=f'(t)-h'(t)f'(h(t))=0
\label{F1_2}
\end{eqnarray}
(\ref{Hypo_lemme_4}) implies $f'$ is a non-zero real analytic function whose zero are isolated, as well as $f'\circ h$. (\ref{F1_2}) and (\ref{F1_1}) implies $f'\circ h=(f'\circ h)h'$; thus $h'$ is equal to $1$ on the dense subset of $U$ where $f'\circ h\neq0$ and by continuity on the whole open subset $U$. As $U$ is connected, it yields for $h$:
$$\forall t\in U,\ h(t)=t'_{m,2}+(t-t'_{m,1})$$
The following combination of $f$ and $f''$ is a polynomial for the variable $t$:
\begin{eqnarray}
f(t)+\epsilon^2f''(t)=(P(t)+\cos(\frac{t}{\epsilon}))+\epsilon^2(P''(t)-\frac{\cos(\frac{t}{\epsilon})}{\epsilon^2})=P(t)+\epsilon^2P''(t)
\label{Poly}
\end{eqnarray}
$t\to((P+\epsilon^2P'')(t)- (P+\epsilon^2P'')(t'_{m,2}+(t-t'_{m,1}))$ is a null polynomial by (\ref{oukoui}) which implies either $P$ to be constant either $t'_{m,2}=t'_{m,1}$; it is in contradiction with the hypothesis of the lemma. \\ 
Finally  $\Sigma_{i=1}^l m_i$ has to be finite and thereby $l$.\\\\

\noindent $\mathbf{Lemma\ 5}$. Let $a<b$, $f\in C^\infty[a,b]$ verifying:
\begin{eqnarray}
\forall t \in [a,b], (f'(t),f''(t))\neq (0,0)
\label{Hypo_lemme_5}
\end{eqnarray}
Suppose there exists $E$ in $\mathbb{R}^3$ such as  the equation for the variable $t$, $E=f^{*3}(t)$ admits $m\geqslant 2$ solutions distinct $(t_i)_{i\leqslant m}$ with $m$ finite. (\ref{Hypo_lemme_5}) ensures $E\neq0$.\\
Suppose there exists $\beta>0$, such as $f$ is real-analytic on $ \sqcup_{j=1}^m ]t_j-\beta,t_j+\beta[$.\\
For any $j>1$, (\ref{Hypo_lemme_5}) ensures $\|f^{1*4}(t_j)\|\neq 0$; thus the sequence $(\alpha_{j,n})_n$ is well-defined by induction on $n$:
\begin{eqnarray}
\begin{array}{lll}
\alpha_{j,1}&=&\frac{f^{1*4}(t_j)\cdot  f^{1*4}(t_1)}{\|f^{1*4}(t_j)\|^2}\\\\
\alpha_{j,n}&=&n!\frac{f^{1*4}(t_j)\cdot \Big(\frac{f^{n*(n+3)}(t_1)}{n!}-\sum\limits_{\sum_{k=1}^{n-1} k m_k=n}\Big[\prod_{l=1}^{n-1} \frac{\alpha_{j,l}^{m_l}}{(m_l! l!^{m_l}) } \Big]f^{(\sum_{k=1}^{n-1} m_k)*(\sum_{k=1}^{n-1} m_k+3)}(t_j) \Big)}{\|f^{1*4}(t_j)\|^2}
\end{array}
\label{Cons_alpha}
\end{eqnarray}
as well as the sequence of planes or straight lines:
\begin{eqnarray}
P_{j,n}=\mathbb{R}f^{1*4}(t_j)+\mathbb{R}\Big(\frac{f^{n*(n+3)}(t_1)}{n!}-\sum  \limits_{\sum_{k=1}^{n-1} k m_k=n}\Big[\prod_{l=1}^{n-1} \frac{\alpha_{j,l}^{m_l}}{(m_l! l!^{m_l}) } \Big] f^{(\sum_{k=1}^{n-1} m_k)*(\sum_{k=1}^{n-1} m_k+3)}(t_j)\Big)
\label{Cons_Pjn}
\end{eqnarray}
\\
Let 
\begin{eqnarray}
V_0\in\mathbb{R}^3\setminus \bigcup\limits_{n=1}^\infty \Big(\bigcup \limits_{j=1}^m (P_{j,n}) \Big).
 \label{Hyp1_lemme_5}
\end{eqnarray}
Let $\tau_1$, $\tau_2$, $\delta$ three of the variable $\epsilon_1$ real analytic functions on a neighborhood of 0 such as 
\begin{eqnarray}
 \tau_1(0)=0,\; \tau_2(0)=0,\; \delta^{*3}(t_1)=V_0.
 \label{Hyp2_lemme_5}
 \end{eqnarray}
The coefficients of the analytical expansions of $\tau_1$, $\tau_2$ in $0$ are noted $(\tau_{1,n})$, $(\tau_{2,n})$.
Let $j\in\{2,\cdots,m\}$ and suppose: 
\begin{eqnarray}
\forall \epsilon_1\in ]-\beta,\beta[, f^{*3}(t_j+\tau_2(\epsilon_1))-f^{*3}(t_1+\epsilon_1)=\tau_1(\epsilon_1) \delta^{*3}(t_1+\epsilon_1)
\label{Hyp3_lemme_5}
\end{eqnarray}
Then for any $n\in\mathbb{N}^*,$ $P_{j,n}$ is a straight line and
\begin{eqnarray}
&&\tau_{1,n}=0 \label{Hypi1_lemme5}\\
&&\tau_{2,n}=\alpha_{j,n}\label{Hypi2_lemme5}\label{Hypi2_lemme5}
\end{eqnarray}

\noindent $\mathbf{Proof\ of\ Lemma\ 5}$. Let $j\in\{2,\cdots,m\}$.  Let $\epsilon_1\in ]-\beta,\beta[$. Let expand the left side of (\ref{Hyp3_lemme_5}) in $\epsilon_1$:
\begin{eqnarray*}
f^{*3}(t_j+\tau_2(\epsilon_1))-f^{*3}(t_1+\epsilon_1)=\sum \limits_{n=1}^{\infty}\Big(\frac{(f(t_j+\tau_2))^{n*(n+3)}(0)}{n!} \epsilon_1^n -\frac{f^{n*(n+3)}(t_1)}{n!} \epsilon_1^n\Big)\\
\end{eqnarray*}
Once expanded by the Faa di Bruno formula it gives:
\begin{eqnarray*}
&&f^{*3}(t_j+\tau_2(\epsilon_1))-f^{*3}(t_1+\epsilon_1)=\\
&&\;\;\;\sum \limits_{n=1}^{\infty}\Big(\frac{\sum  \limits_{\sum_{k=1}^{n} k m_k=n}n!\Big[\prod_{l=1}^{n} \frac{\tau_{2,l}^{m_l}}{(m_l! l!^{m_l}) } \Big]f^{(\sum_{k=1}^{n} m_k)*(\sum_{k=1}^{n} m_k+3)}(t_j)}{n!} -\frac{f^{n*(n+3)}(t_1)}{n!} \Big)\epsilon_1^n
\end{eqnarray*}
which can be simplified into:
\begin{eqnarray*}
&&f^{*3}(t_j+\tau_2(\epsilon_1))-f^{*3}(t_1+\epsilon_1)=\\
&&\;\;\;\sum \limits_{n=1}^{\infty}\Big(\sum  \limits_{\sum_{k=1}^{n} k m_k=n}\Big[\prod_{l=1}^{n} \frac{\tau_{2,l}^{m_l}}{(m_l! l!^{m_l}) } \Big] f^{(\sum_{k=1}^{n} m_k)*(\sum_{k=1}^{n} m_k+3)}(t_j) -\frac{f^{n*(n+3)}(t_1)}{n!} \Big)\epsilon_1^n
\end{eqnarray*}
Let the sequence $(V_k)_{k\in \mathbb{N}}$ the analytic coefficient of $\delta^{*3}$. Let expand the right side of (\ref{Hyp3_lemme_5}):
\begin{eqnarray*}
\tau_1(\epsilon_1) \delta^{*3}(t_1+\epsilon_1)=\sum \limits_{n=1}^{\infty}\Big(\sum \limits_{k=0}^{n}\tau_{1,n-k}V_k \Big)\epsilon_1^n
\end{eqnarray*}
By identification of the coefficients of the analytic expansion on both sides of (\ref{Hyp3_lemme_5}):
\begin{eqnarray}
\forall n\in\mathbb{N}^*, \sum \limits_{k=0}^{n}\tau_{1,n-k}V_k=\sum  \limits_{\sum_{k=1}^{n} k m_k=n}\Big[\prod_{l=1}^{n} \frac{\tau_{2,l}^{m_l}}{(m_l! l!^{m_l}) } \Big] f^{(\sum_{k=1}^{n} m_k)*(\sum_{k=1}^{n} m_k+3)}(t_j) -\frac{f^{n*(n+3)}(t_1)}{n!} \label{identif}
\end{eqnarray}
For $n=1$, (\ref{identif}) reads:
\begin{eqnarray*}
\tau_{1,1}V_0=\tau_{2,1}f^{1*4}(t_j) -f^{1*4}(t_1)
\end{eqnarray*}
As $V_0$ lays outside of $P_{j,1}$, (\ref{Hypi1_lemme5}, \ref{Hypi2_lemme5}) are verified for $n=0$:
\begin{eqnarray*}
\tau_{1,1}&=&0\\ 
\tau_{2,1}&=&\alpha_{j,1}
\end{eqnarray*}
and $P_{j,1}$ is a straight line.
Let suppose the induction hypothesis (\ref{Hypi1_lemme5}, \ref{Hypi2_lemme5}) have been proved until the order $n-1$. We have:
\begin{eqnarray*}
\sum \limits_{k=0}^{n}\tau_{1,n-k}V_k&=&\frac{\partial^n (f^{*3}(t_j+\tau_2))}{n!\partial \epsilon_1^n} \Big|_{\epsilon_1=0}-\frac{f^{n*(n+3)}(t_1)}{n!}
\end{eqnarray*}
By hypothesis of induction, for $k>1$, $\tau_{1,n-k}=0$, (\ref{identif}) rewrites:
\begin{eqnarray*}
\tau_{1,n}V_0&=&\frac{f^{1*4}(t_j)}{n!}\tau_{2,n}+\sum  \limits_{\sum_{k=1}^{n-1} k m_k=n}\Big[\prod_{l=1}^{n-1} \frac{\tau_{2,l}^{m_l}}{(m_l! l!^{m_l}) } \Big] f^{(\sum_{k=1}^{n-1} m_k)*(\sum_{k=1}^{n-1} m_k+3)}(t_j) -\frac{f^{n*(n+3)}(t_1)}{n!}
\end{eqnarray*}
By hypothesis of induction, for $k<n$, $\tau_{2,k}=\alpha_{j,k}$
\begin{eqnarray}
\tau_{1,n}V_0&=&\frac{f^{1*4}(t_j)}{n!}\tau_{2,n}+\sum  \limits_{\sum_{k=1}^{n-1} k m_k=n}\Big[\prod_{l=1}^{n-1} \frac{\alpha_{j,l}^{m_l}}{(m_l! l!^{m_l}) } \Big]f^{(\sum_{k=1}^{n-1} m_k)*(\sum_{k=1}^{n-1} m_k+3)}(t_j) -\frac{f^{n*(n+3)}(t_1)}{n!}\label{cucu}
\end{eqnarray}
As $V_0$ lays outside of $P_{j,n}$
\begin{eqnarray*}
\tau_{1,n}&=&0\\
\tau_{2,n}&=&\alpha_{j,n}\\
\end{eqnarray*}
Moreover because of (\ref{cucu}), $P_{j,n}$ is a straight line.\\  
Both the induction step and the initial step stand; the proof is always true.\\\\

\noindent $\mathbf{Lemma\ 6.}$ Let $a<b$, $f \in C^\infty([a,b])$ verifying the four following hypothesis:
\begin{eqnarray}
\forall t \in [a,b], (f'(t),f''(t))\neq (0,0)
\label{Hypo1_lemme_6}
\end{eqnarray}
There exists $E$ in $\mathbb{R}^3$ such as  the equation for the variable $t$, $E=f^{*3}(t)$ admits $m\geqslant 2$ solutions distinct $(t_i)_{1\leqslant i\leqslant m}$ with $m$ finite. (\ref{Hypo1_lemme_6}) ensures $E\neq0$.\\
The total number $M$ of  couples $(t'_{i,1},t'_{i,2})_{i\leqslant M}$ such as $f^{*3}(t'_{i,1})=f^{*3}(t'_{i,2})$ with  $t'_{i,1}\neq t'_{i,2}$ in $[a,b]$ is finite.\\
There exists $\beta_{max}>0$, such as $f$ is real-analytic on $\sqcup_{j=1}^m ]t_j-\beta_{max},t_j+\beta_{max}[$. \\
Then there exist $\rho_0>0$ and $\beta_{0}>0$ such as for any $\rho\in]0,\rho_0]$, for any $\beta\in]0,\beta_{0}]$, a linear combination $\delta$ of plateau functions real-analytic on their open support can be constructed such that:
\begin{eqnarray}
Supp(\delta) =   ]t_1-\beta,t_1+\beta[
\label{conc1_Lemme6}
\end{eqnarray}
\begin{eqnarray}
\|\delta\|_s<\rho
\label{conc2_Lemme6}
\end{eqnarray}
\begin{eqnarray}
\forall  (t'_{1},t'_{2})  \in  [t_1-\beta,t_1+\beta]\times [a,b],\  t'_{1}\neq t'_{2}\Rightarrow(f+\delta)^{*3}(t'_1)\neq  (f+\delta)^{*3}(t'_2)
\label{conc3_Lemme6}
\end{eqnarray}

\noindent $\mathbf{Proof\ of\ Lemma\ 6.}$ As $M$ is finite, there exists  $\rho_0>0$ such as the ball $\mathbf{B}(f^{*3}(t_1)),\rho_0)$ contains no other self-intersection of $f^{*3}$. The number of connected components of $f^{*3}-\{f^{*3}(t_1)\}$ included in $\mathbf{B}(f^{*3}(t_1)),\rho_0)$ is exactly $m$. Let $\rho\in ]0,\rho_0[$. Let $\beta_0 \in ]0,\beta_{max}]$ such as for any $i\leqslant m$, $f^{*3}$ remains in  $\mathbf{B}(f^{*3}(t_i),\rho)$ on $[t_i-\beta_0,t_i+\beta_0]$. Let $\beta\in ]0,\beta_0[$. \\
For any $j>1$, let define $\alpha_{j,n},\; P_{j,n}$ by induction on $n$ according to the definition (\ref{Cons_alpha}), (\ref{Cons_Pjn}). Let 
\begin{eqnarray}
V_0\in\mathbb{R}^3\setminus \bigcup\limits_{n=1}^\infty \Big(\bigcup \limits_{j=1}^m (P_{j,n}) \Big).\label{V_00000}
\end{eqnarray}
Let choose three strictly positive real numbers $(\gamma_{j})_{1\leqslant j\leqslant 3}$ such as $\gamma_1\neq\gamma_2$.
\begin{eqnarray*}
\left(
\begin{array}{ccc}
 ^t\psi_{\beta,t,\gamma_1}^{*3}(t)&  ^t\psi_{\beta,t,\gamma_2}^{*3}(t)&  ^t\psi_{\beta,t,\gamma_3}^{1*4}(t)
\end{array}
\right)=\left(
\begin{array}{ccc}
1&1&0\\
0&0&-\gamma_3 \\
-\gamma_1&-\gamma_2&0\\
\end{array}
\right)
\end{eqnarray*}
is invertible. 
We can thus note: 
$$v=\left(
\begin{array}{ccc}
  ^t\psi_{\beta,t_1,\gamma_1}^{*3}(t)&  ^t\psi_{\beta,t_1,\gamma_2}^{*3}(t)&  ^t\psi_{\beta,t_1,\gamma_3}^{1*4}(t)
\end{array}
\right)^{-1}V_0.$$
We define $\delta$ for $u \in\ I$:
\begin{eqnarray*}
\delta(u)&=&\sum_{j=1}^2v_j \psi_{\beta,t_1,\gamma_{j}}(u)+v_3\frac{\partial \psi_{\beta,t_1,\gamma_{3}}}{\partial t}(u)
\end{eqnarray*}
by construction 
$$\delta^{*3}(t_1)=V_0.$$ 
As the $(\gamma_{j})_{1\leqslant j\leqslant 3}$ are strictly positive, $\delta$ is $C^\infty$ on $[a,b]$, moreover it is real-analytic on $]t_1-\beta,t_1+\beta[$ and $0$ elsewhere.
For $\epsilon\in]0,\frac{\rho}{\|\delta\|_s}[$ and $u\in [a,b]$, let
\begin{eqnarray*}
g(u,\epsilon)&=&f(u)+\epsilon\delta(u)
\end{eqnarray*}
By choice of $\epsilon$ for any $u\in [a,b]$
\begin{eqnarray}
\|g(u,\epsilon)-f(u)\|\leqslant\ \epsilon\|\delta\|_s<\frac{\rho\|\delta\|_s}{\|\delta\|_s}<\rho_0
\label{bornu}
\end{eqnarray}
\\
\noindent Let suppose there exists an injective sequence $(\epsilon_n)_{n\in \mathbb{N}}$ of  $]0,\frac{\rho}{\|\delta\|_s}[$ whose limit is $0$ and a sequence $(t'_{1,n},t'_{2,n})_{n\in \mathbb{N}}\in [t_1-\beta,t_1+\beta]\times  [a,b]$ such as:
\begin{eqnarray}
t'_{1,n}\neq t'_{2,n},\ g^{*3}(t'_{1,n},\epsilon_n)=g^{*3}(t'_{2,n},\epsilon_n)\label{iuppopopooo}
\label{refoppp}
\end{eqnarray}
As the whole sequence lays in the compact set $[t_1-\beta,t_1+\beta]\times   [a,b]$, there exists a subsequence $(\tilde{t}_{1,n},\tilde{t}_{2,n})$ converging toward $(t'_1,t'_2)\in [t_1-\beta,t_1+\beta]\times  [a,b]$.\\
$(t'_1,t'_2)$ verifies by continuity:
$$f^{*3}(t'_{1})=g^{*3}(t'_{1},0)=g^{*3}(t'_{2},0)=f^{*3}(t'_{2}).$$
The choice of $\beta$ at the beginning of proof implies $t'_{1}=t_1$ and  $t'_{2}=t_j$ for some $j\in\{1,\cdots,m\}$. \\
(\ref{refoppp}) can be rewritten:
\begin{eqnarray}
\forall n\in\mathbb{N},\:f^{*3}(\tilde{t}_{1,n})-f^{*3}(\tilde{t}_{2,n})&=&\epsilon_n(\delta^{*3}(\tilde{t}_{2,n})-\delta^{*3}(\tilde{t}_{1,n}))
\label{refopppuop}
\end{eqnarray}
Let first  suppose $j=1$. A first order development of both sides of the equality gives for any $n\in\mathbb{N}$
\begin{eqnarray*}
f^{*3}(\tilde{t}_{1,n})-f^{*3}(\tilde{t}_{2,n})&=&f^{1*4}(t_{1})(\tilde{t}_{1,n}-\tilde{t}_{2,n})+o(\tilde{t}_{1,n}-\tilde{t}_{2,n})\\
\delta^{*3}(\tilde{t}_{1,n})-\delta^{*3}(\tilde{t}_{2,n})&=&\delta^{1*4}(t_{1})(\tilde{t}_{1,n}-\tilde{t}_{2,n})+o(\tilde{t}_{1,n}-\tilde{t}_{2,n})
\end{eqnarray*}
by identification of the first order term in $(\tilde{t}_{1,n}-\tilde{t}_{2,n})$ in (\ref{refopppuop}), $f^{1*4}(t_{1})=-\epsilon_n \delta^{1*4}(t_{1})$. As the limit of $\epsilon_n $ is zero, $f^{1*4}(t_{1})=0$ which is in contradiction with (\ref{Hypo1_lemme_6}).\\\\
Let suppose $j>1$. Provided $n$ is sufficiently high, $\tilde{t}_{2,n}$ is close enough of $t_j$ to be outside the support of $\delta$, so $\delta^{*3}(\tilde{t}_{2,n})=0$. \\
Let $$V_{0,\bot}=f^{1*4}(t_{j})-\frac{f^{1*4}(t_{j})\cdot V_0}{(V_0)^2}V_0.$$
By construction,  $V_{0,\bot}$ is orthogonal to $V_0$. Moreover $V_0$ has been chosen linearly independent of $f^{1*4}(t_{j})$ (\ref{V_00000}), thus the following dot product is non-zero:
\begin{eqnarray}
V_{0,\bot}\cdot f^{1*4}(t_{j})=f^{1*4}(t_{j})^2-\frac{(f^{1*4}(t_{j})\cdot V_0)^2}{(V_0)^2}\neq 0\label{doto}
 \end{eqnarray}
We note $\widetilde{\beta}\in]0,\beta]$ such as $\delta^{*3}(]t_1-\widetilde{\beta},t_1+\widetilde{\beta}[)\cdot V_0$ does not contain $0$. By construction $\delta$ is not real-analytic on $[a,b]$ but is real-analytic on $]t_1-\beta,t_1+\beta[$ and thus on $]-\widetilde{\beta},\widetilde{\beta}[$. We define the following real-analytic functions $F$, $G$, $H$:
\begin{eqnarray*}
F&:&\Big\{ \begin{array}{ccc}]-\widetilde{\beta},\widetilde{\beta}[^2&\to&\mathbb{R}\\
(\varepsilon_1,\varepsilon_2)&\to&\frac{(f^{*3}(t_j+\varepsilon_2)-f^{*3}(t_1+\varepsilon_1))\cdot V_0}{ \delta^{*3}(t_1+\varepsilon_1)\cdot V_0}
 \end{array}\\
G&:&\Big\{ \begin{array}{ccc}]-\widetilde{\beta},\widetilde{\beta}[^2&\to&\mathbb{R}^3\\
(\varepsilon_1,\varepsilon_2)&\to&(f^{*3}(t_1+\varepsilon_1)-f^{*3}(t_j+\varepsilon_2)-F(\varepsilon_1,\varepsilon_2) \delta^{*3}(t_1+\varepsilon_1))
 \end{array}\\
H&:&\Big\{ \begin{array}{ccc}]-\widetilde{\beta},\widetilde{\beta}[^2&\to&\mathbb{R}\\
(\varepsilon_1,\varepsilon_2)&\to& G(\varepsilon_1,\varepsilon_2)\cdot V_{0,\bot}
 \end{array}
\end{eqnarray*}
By construction of $F$, $-\epsilon_{n}=F(\tilde{t}_{1,n}-t_{1},\tilde{t}_{2,n}-t_{j})$.\\
As $f^{*3}(t_1)=f^{*3}(t_j)$, $G(0,0)=0$ and $H(0,0)=0$. The $\varepsilon_2$-derivative of $H$
\begin{eqnarray*}
\partial_{\varepsilon_2}H(\varepsilon_1,\varepsilon_2)=-f^{1*4}(t_j+\varepsilon_2)\cdot V_{0,\bot}+f^{1*4}(t_j+\varepsilon_2)\cdot V_0\frac{\delta^{*3}(t_1+\varepsilon_1)\cdot V_{0,\bot}}{\delta^{*3}(t_1+\varepsilon_1)\cdot V_0}\\
\end{eqnarray*}
is non-zero by (\ref{doto}):
\begin{eqnarray*}
\partial_{\epsilon_2}H(0,0)&=&-f^{1*4}(t_j)\cdot V_{0,\bot}+f^{1*4}(t_j)\cdot V_0 \frac{ V_0\cdot V_{0,\bot}}{\|V_0\|^2}\\
\partial_{\epsilon_2}H(0,0)&=&-f^{1*4}(t_j)\cdot V_{0,\bot}\neq 0
\end{eqnarray*}
As $F$ is a quotient of two real analytic functions on $]-\widetilde{\beta},\widetilde{\beta}[^2$ whose denominator is non-zero, $F$ is real analytic on  $]-\widetilde{\beta},\widetilde{\beta}[^2$; $G$, $H$ are real-analytic as well on $]-\beta,\beta[^2$. The implicit function theorem for real-analytic function (Theorem 2.3.1 of cite{Krantz1992}) ensures there exists a connected neighborhood $U\subset ]-\widetilde{\beta},\widetilde{\beta}[$  of $0$, a connected  neighborhood $V\subset ]-\widetilde{\beta},\widetilde{\beta}[$ of $0$, and a unique real analytic function $\tau_2:U\to V$  such that:
\begin{eqnarray*}
\forall \varepsilon_1 \in\ U, \forall \varepsilon_2 \in\ V, \ H(\varepsilon_1,\varepsilon_2)&=&0\iff \varepsilon_2=\tau_2(\varepsilon_1)
\label{Fimpo}
\end{eqnarray*}
Let $\tau_1$, $G_1$ defined by
\begin{eqnarray*}
\tau_1: \Big\{ \begin{array}{ccc}U&\to&\mathbb{R}\\
\varepsilon_1&\to& F(\varepsilon_1,\tau_2(\varepsilon_1))
 \end{array}\\
G_1:\Big\{ \begin{array}{ccc}U&\to&\mathbb{R}^3\\
\varepsilon_1&\to& G(\varepsilon_1,\tau_2(\varepsilon_1))
 \end{array}
\end{eqnarray*}
\noindent There exists $n_0$ such as for any $n>n_0$, $\tilde{t}_{1,n}-t_{1}\in U$ and $\tilde{t}_{2,n}-t_{2}\in V$. By definition of $\tau_2$, 
\begin{eqnarray}
\tilde{t}_{2,n}-t_{j}=\tau_{2}(\tilde{t}_{1,n}-t_{1}).\label{euy_nic}
\end{eqnarray}
The proof will be separated in two cases: $(\tilde{t}_{1,n})$ is dense around $t_1$, there exists $n_1>0$ such as $(\tilde{t}_{1,n})$ is constant equal $t_1$ for $n>n_1$.\\ 
\textbf{Case 1}. Let first suppose $(\tilde{t}_{1,n})$ is dense around $t_1$. For any $n>n_0$, $G_1(\tilde{t}_{1,n}-t_{1})=0_{\mathbb{R}^3}$; $G_1$ is a real-analytic function of one variable which admits an accumulation of $0$ around $0$ implying $G_1=0_{\mathbb{R}^3}$ on $U$. $G_1=0_{\mathbb{R}^3}$  can be rewritten:
\begin{eqnarray*}
\forall \varepsilon_1\in U, f^{*3}(t_j+\tau_2(\varepsilon_1))-f^{*3}(t_1+\varepsilon_1)=\tau_1(\varepsilon_1) \delta^{*3}(t_1+\varepsilon_1)
\end{eqnarray*}
The hypothesis (\ref{Hyp1_lemme_5}), (\ref{Hyp2_lemme_5}), (\ref{Hyp3_lemme_5}) are thus verified; Lemma 5 ensures: 
\begin{eqnarray*}
\forall n\in \mathbb{N},&&\tau_{1,n}=0\\
&& \tau_{2,n}=\alpha_{j,n}
\end{eqnarray*}
So for $n>n_0$, as $-\epsilon_n=\tau_1(\tilde{t}_{1,n}-t_{1},\tau_{2}(\tilde{t}_{1,n}-t_{1}))$, $\epsilon_n=0$ which is in contradiction with the definition of $(\epsilon_n)_{n\in \mathbb{N}}$. \\\\
It implies there exists $\epsilon_0\in ]0,\frac{\rho}{\|\delta\|_s}[$ such as for any $\epsilon\in]0,\epsilon_0[$:
\begin{eqnarray*}
Supp(\epsilon\delta) =  ]t_1-\beta,t_1+\beta[
\end{eqnarray*}
\begin{eqnarray*}
\|\epsilon\delta\|_s<\rho
\end{eqnarray*}
\begin{eqnarray*}
\forall  (t'_{1},t'_{2})  \in  [t_1-\beta,t_1+\beta]\times [a,b],\  t'_{1}\neq t'_{2}\Rightarrow(f+\epsilon\delta)^{*3}(t'_1)\neq  (f+\epsilon\delta)^{*3}(t'_2)
\end{eqnarray*}
(\ref{conc1_Lemme6}), (\ref{conc2_Lemme6}), (\ref{conc3_Lemme6}) are verified for $\epsilon\delta$.\\\\
\textbf{Case 2}. Let suppose there exists $n_1$ such as for any $n>n_1$, $\tilde{t}_{1,n}=t_1$. Because of (\ref{euy_nic}), $\tilde{t}_{2,n}-t_j=\tau_2(0)=0$. By hypothesis $f^{*3}(t_1)=f^{*3}(t_j)$, by construction $\delta^{*3}(t_1)=V_0$ and $\delta^{*3}(t_j)=0$,  (\ref{refopppuop}) rewrites:
\begin{eqnarray*}
\forall n\in\mathbb{N},\: 0=f^{*3}(t_1)-f^{*3}('_j)=-\epsilon_n V_0
\end{eqnarray*}
which is impossible as both $\epsilon_n$ and $V_0$ are non-zero.\\\\

\noindent $\mathbf{Lemma\ 7.}$  Let $a<b$, $\epsilon>0$,  $P$ be a non-constant polynomial on $[a,b]$ and $$f=P+c_\epsilon.$$
We suppose:
\begin{eqnarray}
\forall t \in [a,b], (f'(t),f''(t))\neq (0,0)
\label{Hypo_lemme_7}
\end{eqnarray}
Suppose there exists $E$ in $\mathbb{R}^3$ such as  the equation for the variable $t$, $E=f^{*3}(t)$ admits $m\geqslant 2$ solutions distinct $(t_i)_{1\leqslant i\leqslant m}$. (\ref{Hypo_lemme_7}) ensures $E\neq0$. Let
 \begin{eqnarray*}
\beta_{dis}&=&\frac{\min_{1\leqslant i < k\leqslant m}(|t_{i}-t_{k}|)}{3}\\
\end{eqnarray*}\\
Then for any $\rho>0$, for any $\beta\in ]0,\beta_{dis}[$, there exists a linear combination $\delta$ of plateau functions real-analytic on the interior of their support such that:
\begin{eqnarray}
Supp(\delta) \subset    \sqcup_{j=1}^m ]t_j-\beta,t_j+\beta[
\label{conc1_Lemme8}
\end{eqnarray}
\begin{eqnarray}
\|\delta\|_s<\rho
\label{conc2_Lemme8}
\end{eqnarray}
\begin{eqnarray}
\forall  (t'_{1},t'_{2})  \in   \sqcup_{j=1}^m [t_j-\beta,t_j+\beta]\times [a,b],\  t'_{1}\neq t'_{2}\Rightarrow(f+\delta)^{*3}(t'_1)\neq  (f+\delta)^{*3}(t'_2)
\label{conc3_Lemme8}
\end{eqnarray}
\noindent $\mathbf{Proof\ of\ Lemma\ 7.}$ 
As $f$ is the sum of  a non-constant polynomial $P$ and $c_\epsilon$ and $f$ verifies (\ref{Hypo_lemme_7}), Lemma 4 can be applied. Lemma 4 ensures both that the number $l$ of vectors $(E_{i})_{i\in\{1,\cdots,l\}}\in\ \mathbb{R}^3$ such that the equation
$$E_{i}= f^{*3}(t)$$ 
admits more than one solution is finite and that the number $M$ of  couples $(t'_{i,1},t'_{i,2})_{i\leqslant M}$ such as $f^{*3}(t'_{i,1})=f^{*3}(t'_{i,2})$ with  $t'_{i,1}\neq t'_{i,2}$ in $[a,b]$ is finite. Moreover $f$ is real-analytic on the whole interval $[a,b]$ as it is the sum of a sinusoidal function and a polynomial. We note $f_1=f$, $f_1$ verifies the four hypothesis of Lemma 6 which ensures there exists $\delta_1$ verifying:
\begin{eqnarray*}
Supp(\delta_1) \subset   ]t_1-\beta,t_1+\beta[
\end{eqnarray*}
\begin{eqnarray*}
\|\delta_1\|_s<\frac{\rho_1}{10}
\end{eqnarray*}
\begin{eqnarray*}
\forall  (t'_{1},t'_{2})  \in [t_1-\beta,t_1+\beta]\times [a,b],\  t'_{1}\neq t'_{2}\Rightarrow(f+\delta_1)^{*3}(t'_1)\neq  (f+\delta_1)^{*3}(t'_2)
\end{eqnarray*} .\\
where $\rho_1=\min(\rho,\|f^{*3}\|_i).$\\
Let $k\in\{1,\cdots,m-2\}$, suppose we have constructed for $j\leqslant k$, $\tilde{\delta}_j$ such as:
\begin{eqnarray*}
Supp(\tilde{\delta}_j) \subset  ]t_j-\beta,t_j+\beta[
\end{eqnarray*}
\begin{eqnarray*}
\|\tilde{\delta}_j\|_s<\frac{\rho_1}{10}
\end{eqnarray*}
\begin{eqnarray}
\forall  (t'_{1},t'_{2})  \in   \sqcup_{j=1}^k [t_j-\beta,t_j+\beta]\times [a,b],\  t'_{1}\neq t'_{2}\Rightarrow(f_k)^{*3}(t'_1)\neq  (f_k)^{*3}(t'_2)
\label{yuoipolo1}
\end{eqnarray}
with $f_k=f+\delta_k$ and $\delta_k=\Sigma_{j=1}^k\tilde{\delta}_j$. As the support of the $(\tilde{\delta}_j)$ are disjoint:\\
$$\|\delta_k\|_s=max_{j\in\{1,\cdots,k\}}\|\tilde{\delta}_j\|_s<\frac{\rho_1}{10}.$$
The hypothesis of Lemma 6 are verified for $f_k$:
\begin{itemize}
\item  $f_k$ is smooth on $[a,b]$. By choice of $\rho_1$ 
$$\|f_k\|_i\geqslant \|f\|_i-\|\delta_k\|_s\geqslant \frac{9\|f\|_i}{10}>0,$$ (\ref{Hypo_lemme_7}) is also verified by $f_k$. 
\item  The equation for the variable $t$, $E=f_k^{*3}(t)$ admits a finite number $m-k\geqslant 2$ of solutions distinct $(t_i)_{i\in\{k+1,\cdots, m\}}$.
\item  The total number of  couples $(t'_{k,i,1},t'_{k,i,2})_{i\leqslant M_k}$ such as $f_k^{*3}(t'_{k,i,1})=f_k^{*3}(t'_{k,i,2})$ with  $t'_{k,i,1}\neq t'_{k,i,2}$ in $[a,b]$ is finite.
\item As $f_k=f$ on $\sqcup_{j=k+1}^m ]t_j-\beta,t_j+\beta[$, $f_k$ is real-analytic on $\sqcup_{j=k+1}^m ]t_j-\beta,t_j+\beta[$.
\end{itemize} 
Lemma 6 can be applied to $f_k$. There exists $\tilde{\delta}_{k+1}$  such as:
\begin{eqnarray*}
Supp(\tilde{\delta}_{k+1}) \subset   ]t_{k+1}-\beta,t_{k+1}+\beta[
\end{eqnarray*}
\begin{eqnarray*}
\|\tilde{\delta}_{k+1}\|_s<\frac{\rho_{1}}{10}
\end{eqnarray*}
\begin{eqnarray}
\forall  (t'_{1},t'_{2})  \in \ [t_{k+1}-\beta,t_{k+1}+\beta]\times [a,b],\  t'_{1}\neq t'_{2}\Rightarrow(f_{k}+\tilde{\delta}_{k+1})^{*3}(t'_1)\neq  (f_{k}+\tilde{\delta}_{k+1})^{*3}(t'_2)\label{yuoipolo2}
\end{eqnarray}
We note:
 $\delta_{k+1}= \delta_{k}+\tilde{\delta}_{k+1}$ and $f_{k+1}=f_{k}+\tilde{\delta}_{k+1}=f+\delta_k.$
 \begin{eqnarray}
Supp(\delta_{k+1}) \subset    \sqcup_{j=1}^{k+1} ]t_j-\beta,t_j+\beta[.
\end{eqnarray}
As the $(\tilde{\delta}_{i})$ have disjoint supports
\begin{eqnarray*}
\|\delta_{k+1}\|_s=\max_{i\in\{1,\cdots,k+1\}}\|\tilde{\delta}_{i}\|_s<\rho
\end{eqnarray*}
Let prove (\ref{yuoipolo1}) for $f_{k+1}$ by enumerating 3 cases:\\
- Case 1: For $(t'_{1},t'_{2})  \in   \sqcup_{j=1}^k [t_j-\beta,t_j+\beta]\times [a,b]\setminus[t_{k+1}-\beta,t_{k+1}+\beta],\  t'_{1}\neq t'_{2},$
\begin{eqnarray*}
f_{k+1}^{*3}(t'_1)&=& f_k^{*3}(t'_1)\\
f_{k+1}^{*3}(t'_2)&= &f_k^{*3}(t'_2)
\end{eqnarray*}
By (\ref{yuoipolo1}), $f_{k+1}^{*3}(t'_1)\neq f_{k+1}^{*3}(t'_2)$.\\
- Case 2: For $(t'_{1},t'_{2})  \in  [t_{k+1}-\beta,t_{k+1}+\beta]\times [a,b],\  t'_{1}\neq t'_{2},$
\begin{eqnarray*}
f_{k+1}^{*3}(t'_1)&=&(f_k+\tilde{\delta}_{k+1})^{*3}(t'_1)\\
f_{k+1}^{*3}(t'_2)&= &(f_k+\tilde{\delta}_{k+1})^{*3}(t'_2)
\end{eqnarray*}
By (\ref{yuoipolo2}), $f_{k+1}^{*3}(t'_1)\neq f_{k+1}^{*3}(t'_2)$.\\
- Case 3: For $(t'_{1},t'_{2})  \in   \sqcup_{j=1}^k [t_j-\beta,t_j+\beta]\times[t_{k+1}-\beta,t_{k+1}+\beta],\  t'_{1}\neq t'_{2},$ is proved  in the same way than case 2 by exchanging the role of $t'_1$ and $t'_2$.\\
By induction, a function $\delta_{m-1}$ satisfying (\ref{conc1_Lemme8}, \ref{conc2_Lemme8}, \ref{conc3_Lemme8}) can be constructed.\\\\

\noindent $\mathbf{Lemma\ 8.}$ Let $a<b$, $\epsilon>0$,  $P$ be a non-constant polynomial on $[a,b]$ and $$f=P+c_\epsilon.$$
We suppose:
\begin{eqnarray}
\forall t \in [a,b], (f'(t),f''(t))\neq (0,0)
\label{Hypo_lemme_9}
\end{eqnarray}
Then for any $\rho>0$, there exists a linear combination $\delta$ of plateau functions real-analytic on the interior of their support such as:
\begin{eqnarray}
\|\delta\|_s<\rho
\label{conc2_Lemme9}
\end{eqnarray}
and such  that $g=f+\delta$ verifies:
\begin{eqnarray}
\forall t \in [a,b], (g^{\prime}(t),g^{\prime\prime}(t))\neq (0,0)
\label{conc3_lemme_9}
\end{eqnarray}
\begin{eqnarray}
\forall  (t'_{1},t'_{2})  \in   [a,b]^2,\  t'_{1}\neq t'_{2}\Rightarrow g^{*3}(t'_1)\neq  g^{*3}(t'_2)
\label{conc4_Lemme9}
\end{eqnarray}

\noindent $\mathbf{Proof\ of\ Lemma\ 8.}$
As (\ref{Hypo_lemme_9}) is verified by $f$, Lemma 4 ensures there exists a finite number $l$ of vectors $(E_{i})_{i\in\{1,\cdots,l\}}\in\ \mathbb{R}^3$, such that the equation
$$E_{i}= f_n^{*3}(t)$$ 
admits more than one solution. For each $i\leqslant l$,  Lemma 4 also ensures  there exists $m_{i}$ distinct solutions $(t_{i,j})_{1\leqslant j\leqslant m_{i}}$ with $m_{i}$ finite.\\
\noindent Let introduce some notations :\\
- $\tau$ is defined by: 
\begin{eqnarray*}
\tau&=&\frac{\min_{1\leqslant i \leqslant k\leqslant l,\  j\in\{1,\cdots,m_{i}\},\  l\in\{1,\cdots,m_{k}\}, (i,j)\neq (k,l)}(|t_{i,j}-t_{k,l}|)}{10}\\
\end{eqnarray*}
- $\mu$ is defined by: 
\begin{eqnarray*}
\mu&=&\min_{1\leqslant i< k\leqslant l,\  j\in\{1,\cdots,m_{i}\},\  l\in\{1,\cdots,m_{k}\}}\mu_{i,j,k,l}
\end{eqnarray*}
where for any $ (i,k)\in \{1,\cdots, l\}^2,\  j\in\{1,\cdots,m_{i}\}, l\in\{1,\cdots,m_{k}\},$ we note:
\begin{eqnarray*}
\mu_{i,j,k,l}&=&\inf_{t'_1\in[t_{i,j}-\tau,t_{i,j}+\tau],\ t'_2\in[t_{k,l}-\tau,t_{k,l}+\tau]}dist(f^{*3}(t'_1)-f^{*3}(t'_2))
\end{eqnarray*}
- $\nu=\inf_{[a,b]}\sqrt{f^{\prime2}+f^{\prime\prime2}}$.\\
By construction, for $i=k$,  $\mu_{i,j,i,l}=0$, and for $i\neq k$, $\mu_{i,j,k,l}>0$, thus $\mu$ is strictly positive; by (\ref{Hypo_lemme_9}), $\nu$ is also strictly positive.\\

\noindent $f$ is real-analytic on $[a,b]$ and $l$ is finite as well as the different $(m_i)_{i\in\{1,\cdots,l\}}$, so Lemma 7 ensures for each  $i\leqslant l$, there exists  $\delta_{i}\in C^{\infty}[a,b]$ which verifies the three following properties: \\
- The inclusion: $$Supp(\delta_{i})\subset Supp_{i}=\sqcup_{1\leqslant j\leqslant m_{i}}(]t_{i,j}-\tau,t_{i,j}+\tau[)$$
- The majoration: $$\|\delta_{i}\|_s<\frac{min(\mu,\nu,\rho)}{10}$$
-  For any $(t'_{1},t'_{2})  \in Supp_{i}\times [a,b]$
\begin{eqnarray} 
 t'_{1}\neq t'_{2}\Rightarrow  (f+\delta_{i})^{*3}(t'_1)\neq   (f+\delta_{i})^{*3}(t'_2)\label{Taratata}
 \end{eqnarray}\\
Let define $\delta=\sum_{1\leqslant i \leqslant l} \delta_{i}$, $Supp=\sqcup_{1\leqslant i\leqslant l}Supp_{i}$, and $g=f+\delta$.\\ 
The support of the $(\delta_i)$ are disjoint so:
$$\|\delta\|_s=\max_{i\in\{1,...,l\}}(\|\delta_{i}\|_s)<\rho$$
which means $\delta$ verifies (\ref{conc2_Lemme9}).\\
For $t\in[a,b]$,  (\ref{conc3_lemme_9}) is verified: $$\sqrt{g'(t)^2+g''(t)^2} \geqslant \sqrt{f'(t)^2+f''(t)^2}- \|\delta\|_s\geqslant \frac{9\nu}{10}> 0.$$
Let look if $g$ respects (\ref{conc4_Lemme9}) by enumerating the four possible cases for $(t'_1,t'_2)$: \\
- If  $(t'_1,t'_2)\in \Big([a,b]\setminus Supp\Big)^2$: \\
 $ g^{*3}(t'_1)= f^{*3}(t'_1)$, and $g^{*3}(t'_2)= f^{*3}(t'_2)$. \\
As $f^{*3}$ has no self-intersection on $[a,b]\setminus Supp$, so if $t'_1\neq t'_2$,  then $g^{*3}(t'_1)\neq  g^{*3}(t'_2)$.\\\\
- If  $(t'_1,t'_2)\in\ \Big( Supp_{i}\Big)\times  \Big([a,b]\setminus Supp \Big)$ for $i\in\{1,\cdots,l\}$:\\
 $ g^{*3}(t'_1)= (f+\delta_{i})^{*3}(t'_1)$, and $g^{*3}(t'_2)= f^{*3}(t'_2)=(f+\delta_{i})^{*3}(t'_2)$. \\
By (\ref{Taratata}), $ g^{*3}(t'_1)\neq g^{*3}(t'_2)$.\\\\
- If $(t'_1,t'_2)\in\ \Big( Supp_{i}\Big)^2$ for $i\in\{1,\cdots,l\}$:\\
 $ g^{*3}(t'_1)= (f+\delta_{i})^{*3}(t'_1)$, and $g^{*3}(t'_2)=(f+\delta_{i})^{*3}(t'_2)$. \\
 By (\ref{Taratata}), $ g^{*3}(t'_1)\neq g^{*3}(t'_2)$.\\\\
- If $(t'_1,t'_2)\in\ \Big( Supp_{i}\Big)\times \Big( Supp_{j}\Big)$ for $i\neq j$ in $\{1,\cdots,l\}$:\\
 $ g^{*3}(t'_1)= (f+\delta_{i})^{*3}(t'_1)$, and $g^{*3}(t'_2)=(f+\delta_{j})^{*3}(t'_2)$. \\
 The triangular inequalities give:
\begin{eqnarray*}
 dist(g^{*3}(t'_1),g^{*3}(t'_2))&\geqslant& dist(f^{*3}(t'_1),f^{*3}(t'_2))-\|\delta_{i}\|_s-\|\delta_{j}\|_s\\
 dist(g^{*3}(t'_1),g^{*3}(t'_2))&\geqslant &\mu-\frac{2\mu}{10}=\frac{8\mu}{10}>0
 \end{eqnarray*}
 So $g^{*3}(t'_1)\neq g^{*3}(t'_2)$. \\
 In any of the four cases, $g=f+\delta$ respects (\ref{conc4_Lemme9}). \\

\noindent $\mathbf{Reference}$


\begin{thebibliography}{06}
\bibitem[Rubel(1981)]{Rubel1981} L.A. Rubel, A universal differential equation. Bull. Amer. Math. Soc. (1981),345-349.
\bibitem[Pour-El(1974)]{Pour-El1974} M. B. Pour-El, Abstract computability and its relation to the general purpose analog computer some connections between logic, differential equations, and analog computers). Trans. Amer. Math. Soc. (1974), 199 ,1-28.
\bibitem[Shinar, Feinberg(2010)]{Shinar2010} G. Shinar, M. Feinberg , Structural sources of robustness in biochemical reaction networks. Science. (2010), 327,1389-91.
\bibitem[Boshernitzan(1986)]{Boshernitzan1986} M. Boshernitzan, Universal Formulae and Universal Differential Equations. Ann. Math. (1986),124,273-291
\bibitem[Elsner(1992)]{Elsner1992} C. Elsner, On the approximation of continuous functions by C°-solutions of third-order algebraic differential equations. Math. Nachr., 157 (1992), pp. 235Ð241
\bibitem[Berger, Gostiaux(1992)]{Berger1992} M. Berger, B. Gostiaux, G\'eometrie diff\'erentielle : vari\' et\' es, courbes et surfaces. Presses universitaires de France (1992)
\bibitem[Krantz, Parks 1992]{Krantz1992} S.G. Krantz, H. R. Parks,  A Primer of Real Analytic functions.  Birkhauser (1992)
\end{thebibliography}
\end{document}